\documentclass[a4paper,12pt,reqno]{amsart}
\usepackage{amsmath,amsfonts}
\newtheorem{thm}{Theorem}
\newtheorem{conj}{Conjecture}

\newcommand{\Mbar}{{\overline{\mathcal{M}}}}

\newcommand{\field}[1]{\ensuremath{\mathbb{#1}}}
\newcommand{\del}{\partial}

\begin{document}

\title[Large genus asymptotics of Weil-Petersson volumes]
{On the large genus asymptotics of Weil-Petersson volumes}

\author{Peter Zograf}
\address{Steklov Mathematical Institute,
St. Petersburg 191023 Russia}
\email{zograf@pdmi.ras.ru}

\date{\today}
\thanks{Partially supported by the RFBR grant 08-01-00379-a and 
by the President of Russian Federation grant NSh-2460.2008.1.}

\begin{abstract}
A relatively fast algorithm for evaluating Weil-Peters-son volumes of moduli spaces 
of complex algebraic curves is proposed. On the basis of numerical data, a conjectural 
large genus asymptotics of the Weil-Petersson volumes is computed. 
Asymptotic formulas for the intersection numbers involving
$\psi$-classes are conjectured as well. The accuracy of the formulas is high enough
to believe that they are exact. 
\end{abstract}

\maketitle

The aim of this note is to report on the recent progress in computing Weil-Petersson volumes
of moduli spaces of complex algebraic curves (with or without marked points) that resulted from better programming, software and hardware as compared to \cite{Z}. 
The numerical evidence led us to a plausible guess about their large genus asymptotic behavior
that may have further applications in algebraic geometry, combinatorics, dynamical systems and
string theory.

\section{Algorithms}

Let $\Mbar_{g,n}$ denote the moduli space of stable $n$-pointed genus $g$ complex algebraic 
curves. The universal curve $p\!:\Mbar_{g,n+1}\rightarrow \Mbar_{g,n}$ has $n$ canonical 
sections $x_1,\dots,x_n$ given by the marked points. 
Put  $\psi_i = c_1(x_i^*\omega)\in H^2 (\Mbar_{g,n},\field{Q})$, where $\omega$
is the relative dualizing sheaf on $\Mbar_{g,n+1}$. The first Mumford class of $\Mbar_{g,n}$ is
the direct image class
$$\kappa_1=p_*\psi_{n+1}^2=\int_{\mbox{\tiny{fiber}}}\psi_{n+1}^2\,\in H^2(\Mbar_{g,n},\field{Q}).$$

The Weil-Petersson metric is K\"{a}hler on $\mathcal{M}_{g,n}$. Its symplectic form $\omega_{WP}$
extends to $\Mbar_{g,n}$ as a closed current and represents the class 
$2\pi^2\kappa_1\in H^2(\Mbar_{g,n},\field{R})$ (see \cite{Wo}). By definition, the (normalized) 
Weil-Petersson volume of  $\mathcal{M}_{g,n}$ is just its standard symplectic volume
with respect to the form $\frac{1}{2\pi^2}\omega_{WP}$:
\begin{equation}
V_{g,n}=\frac{1}{(3g-3+n)!}\int_{\Mbar_{g,n}}\kappa_1^{3g-3+n}.
\end{equation}
For all $g,n\geq 0$ with $2g+n\geq 3$ these are positive rational numbers.

Below we describe an algorithm for computing Weil-Petersson volumes (see \cite{Z} for details).

\begin{thm}
Let 
$$\del_0=\frac{1}{t}\left(\frac{\del}{\del y}
-\frac{x(y)}{y}\frac{\del}{\del t} \right), \quad\quad
\del_1=-\frac{\del}{\del t}+y\del_0,$$
where
$$x(y)=-\sqrt{y}J'_0(2\sqrt{y})=\sum_{k=1}^\infty\frac{(-1)^{k-1}}{(k-1)!}\;\frac{y^k}{k!}$$
($J_0$ denotes the Bessel function of the first kind).\\ 
Then\\
(i) the KdV equation
\begin{equation*}
\del_1 u=\del_0\left(\frac{u^2}{2}+h^2\frac{\del_0^2 u}{12}\right)
\end{equation*}
has a unique solution  $u(y,t)=y+h^{2g}\sum_{g=1}^\infty u_g(y,t)$ where
each $u_g(y,t)$ is a Laurent polynomial in $t$ of the form
$$u_g(y,t)=\sum_{k=2g+1}^{5g-1}u_{g;k}(y)t^{-k};$$\\
(ii) for each $g\geq 2$ the equation 
$$\del_0^2\phi_g(y,t)=u_g(y,t)$$ 
has a unique solution of the form 
$$\phi_g(y,t)=\sum_{k=2g-1}^{5g-5}\phi_{g;k}(y)t^{-k};$$
(iii) for any $g,n\geq 0$ the Weil-Petersson volume of $\mathcal{M}_{g,n}$ is given by the formula
$$V_{g,n}=\del_0^n\phi_g(y,t)\left|_{y=0,\,t=1}\right. .$$
\end{thm}

The above theorem extends to the intersection numbers involving $\psi$-classes.
Fix a set $d=(d_1,\dots,d_n)$ of non-negative integers and put $|d|=d_1+\dots +d_n.$ 
Consider 
\begin{equation}
V_{g,n;d}=\frac{1}{(3g-3+n-|d|)!}
\int_{\Mbar_{g,n}}\psi_1^{d_1}\dots\psi_n^{d_n}\kappa_1^{3g-3+n-|d|}.
\label{psi}\end{equation}

\begin{thm}
Let 
\begin{eqnarray*}
&&\del_0=x_1\left(\frac{\del}{\del y}+x_1^2\left(x_2+\frac{x(y)}{y}
\right)\frac{\del}{\del x_1}+\sum_{k=2}^{\infty}x_{k+1}\frac{\del}{\del x_k}\right), \\
&&\del_1=x_1^2\,\frac{\del}{\del x_1}\,+\,y\,\del_0,
\end{eqnarray*}
with $x(y)=-\sqrt{y}J'_0(2\sqrt{y})$ as above.\\
Then\\
(i) the KdV equation
\begin{equation*}
\del_1 v=\del_0\left(\frac{v^2}{2}+h^2\frac{\del_0^2 v}{12}\right)
\end{equation*}
has a unique solution 
$$v(y,x_1,x_2,\dots)=y+h^{2g}\sum_{g=1}^\infty v_{g}(y,x_1,x_2,\dots),$$ 
where each $v_g(y,t)$ is a polynomial in $x_1$ of the form
$$v_g(y,x_1,x_2,\dots)=\sum_{k=2g+1}^{5g-1}v_{g;k}(y,x_2,x_3,\dots)x_1^{k};$$\\
(ii) for each $g\geq 2$ the equation 
$$\del_0^2\psi_g(y,x_1,x_2,\dots)=v_g(y,x_1,x_2,\dots)$$ 
has a unique solution of the form 
$$\psi_g(y,x_1,x_2,\dots)=\sum_{k=2g-1}^{5g-5}\psi_{g;k}(y,x_2,x_3,\dots)x_1^{k};$$
(iii) the intersection number $V_{g,n;d}$ is given by the formula
$$V_{g,n;d}=\frac{\del^{l_2+l_3+\dots}}{\del x_2^{l_2}\del x_3^{l_3}\dots}\left(x_1^2\frac{\del}{\del x_1}\right)^{l_1}
\del_0^{\,l_0}\psi_g(y,x_1,x_2,\dots)\left|_{y=0,\,x_1=1,\,x_2=x_3=\dots =0}\right. ,$$
where $l_k$ is the number of $d_i$'s equal to $k$.
\end{thm}

The proof follows the same lines as that of Theorem 1 and utilizes an observation of
M.~Kazarian on how to explicitly express mixed intersection numbers of $\psi$- and
$\kappa$-classes in terms of intersection numbers of $\psi$-classes alone \cite{K}.
The details will appear elsewhere.
Note that for $d_1=\dots =d_n=0$ it reduces to Theorem 1 with the obvious change
of variable $x_1=1/t$. 

The main advantage of our algorithm is its speed, and in this
respect it is superior to the algorithms of C.~Faber \cite{F} and M.~Kazarian
\cite{K}, though it loses to both of them in generality.

\section{Asymptotics}

It may be instructive to begin with the large $n$ asymptotics
of Weil-Petersson volumes. The following 
exact asymptotic formula was proven in \cite{MZ} for any fixed $g$: 
\begin{equation}
V_{g,n}=n!\,C^n\, n^{(5g-7)/2}\left(a_g+O\left(1/n\right)\right),
\quad n\rightarrow\infty,
\label{n}\end{equation}
where $C=-z_0J'_0(z_0)$ and $z_0$ is the first positive zero of the Bessel
function $J_0(z)$. The coefficients $a_g$ can also be explicitly computed \cite{MZ}
(in fact, one can even get the complete asymptotic expansion of $V_{g,n}$ as $n\rightarrow\infty$).

The problem seems more challenging when $n$ is fixed and $g\rightarrow\infty$.
We implemented the algorithm of Theorem 1 in a Maple$^{\mbox{\tiny TM}}$
\footnote{\copyright Maplesoft, a division of Waterloo Maple Inc.} program and
computed all numbers $V_{g,n}$ for $g\leq 50$ and $1\leq n\leq 4$. These data led us to

\begin{conj}
For any fixed $n\geq 0$
$$V_{g,n}=(2g)!\left(\frac{2}{\pi^2}\right)^g\,g^{n-7/2}\,\,
\frac{2^{2n-6}}{\sqrt{\pi}}\,\left(1+\frac{c_n}{g}+O\left(\frac{1}{g^2}\right)\right),
\quad g\rightarrow\infty.$$ 
\end{conj}

This formula agrees with the earlier results of \cite{G,ST}. Approximate values of the
constants $c_n$ with $n\leq 4$ are: $c_0 \approx 1.8,\;c_1 \approx 0.75,\; 
c_2 \approx 0.1,\; c_3 \approx -0.15,\; c_4 \approx -0.001.$

Our Maple$^{\mbox{\tiny TM}}$ implementation of the algorithm of Theorem 2 
evaluates the intersection numbers $V_{g,n;d}$ given by (\ref{psi}).
In particular, we computed all $V_{g,n;d}$ with $g\leq 40$ and $l_k\leq 2,\; k=1,2,3$), 
and from that we get

\begin{conj}
For any fixed $n>0$ and a fixed set $d=(d_1,\dots d_n)$ of non-negative integers
$$\lim_{g\rightarrow\infty}\frac{V_{g,n;d}}{V_{g,n}}=\prod_{k\geq 1}\left(\frac{\pi^{2k}}{2^k\,(2k+1)!!}\right)^{l_k},$$ 
where $l_k$ denotes the number of $d_i$'s that are equal to $k$.
\end{conj}

Both these conjectures hold numerically with high accuracy,
so there is a good reason to believe that they are actually true. However, at the moment these asymptotic formulas lack theoretical 
justification.\footnote{As it was proven in M.~Mirzakhani, P.~Zograf, Towards large genus asymptotics of intersection numbers on moduli spaces of curves, Geom. Funct. Anal. {\bf 25}, 1258-1289 (2015), Conjecture 1 holds up to a universal factor that is numerically equal to 1 (cf. Theorem 1.2), 
and Conjecture 2 is true for any set $d=(d_1,\dots d_n)$ (cf. Theorem 4.1).} The next section contains some data and heuristics.

\section{Numerics}

It is known that the order of magnitude of $V_{g,n}$ is $(2g)!$ for large $g$ (see \cite{G,ST}),
so the problem was to find the asymptotics up to the factors of smaller order.
A question of M.~Mirzakhani about the behavior of the ratio
$V_{g-1,n+2}/V_{g,n}$ as $g\rightarrow\infty$ served us as a starting point. 
Computations show that this
ratio decreases with $g$ for $n=0,1$, and increases for any fixed $n\geq 2$. At the 
same time for any fixed $g$ this ratio decreases when $n$ grows. Below is the table
of decimal approximations (rounded up to 10 digits)
for $n=1,2$ and $41\leq g \leq 50$:

\bigskip
\noindent
\begin{center}
\begin{tabular}{|c|c|c|}
\hline
$g$ & \;$V_{g-1,3}/V_{g,1}\;$ & $\;V_{g-1,4}/V_{g,2}$\;\\
\hline \hline
41 & 19.78811999 & 19.68915456 \\
\hline
42 & 19.78695294 & 19.69037678 \\
\hline
43 & 19.78584026 & 19.69154073 \\
\hline
44 & 19.78477824 & 19.69265049 \\
\hline
45 & 19.78376350 & 19.69370974 \\
\hline
46 & 19.78279294 & 19.69472186 \\
\hline
47 & 19.78186376 & 19.69568993 \\
\hline
48 & 19.78097335 & 19.69661676 \\
\hline
49 & 19.78011934 & 19.69750494 \\
\hline
50 & 19.77929954 & 19.69835682 \\
\hline
\end{tabular}
\end{center}
\medskip

The behavior of $V_{g-1,n+2}/V_{g,n}$ suggests that there is a limit as $g\rightarrow\infty$
independent of $n$. D.~Zagier numerically identified this limit with $2\pi^2=19.7392088\dots$ 
(private communication).

The next step is to analyze the behavior of the ratio $2gV_{g,n-1}/V_{g,n}$:

\bigskip
\noindent
\begin{center}
\begin{tabular}{|c|c|c|c|}
\hline
$g$ & \;$2gV_{g,1}/V_{g,2}\;$ & $\;2gV_{g,2}/V_{g,3}\;$ & $\;2gV_{g,3}/V_{g,4}$\;\\
\hline \hline
40 & 0.5082406948 & 0.5031382404 & 0.4981366818 \\
\hline
41 & 0.5080365079 & 0.5030613837 & 0.4981822417 \\
\hline
42 & 0.5078421948 & 0.5029882014 & 0.4982256270 \\
\hline
43 & 0.5076570564 & 0.5029184361 & 0.4982669896 \\
\hline
44 & 0.5074804578 & 0.5028518541 & 0.4983064678 \\
\hline
45 & 0.5073118214 & 0.5027882423 & 0.4983441875 \\
\hline
46 & 0.5071506208 & 0.5027274063 & 0.4983802636 \\
\hline
47 & 0.5069963746 & 0.5026691682 & 0.4984148013 \\
\hline
48 & 0.5068486423 & 0.5026133653 & 0.4984478969 \\
\hline
49 & 0.5067070199 & 0.5025598478 & 0.4984796388 \\
\hline
\end{tabular}
\end{center}
\medskip

It is not hard to see that expected limit in each of the columns is 1/2 as $g\rightarrow\infty$.

These two observations combined together give the factor $2^{g+2n}/\pi^{2g}$ in the 
Weil-Petersson volume asymptotics. Similar to (\ref{n}) it is natural to assume that
the ratio $\frac{\pi^{2g}V_{g,n}}{(2g)!\,2^{g+2n}}$ behaves like $a_n g^{b_n}$, and it
works. Moreover, it appears that $a_n=2^{-6}\pi^{-1/2}$ is independent of 
$n$ and $b_n=n-7/2$. In the table below 
$C_{g,n}=(2g)!\,\left(\frac{2}{\pi^2}\right)^g\,g^{n-7/2}\;\frac{2^{2n-6}}{\sqrt{\pi}}\,$:

\bigskip
\noindent
\begin{center}
\begin{tabular}{|c|c|c|c|c|}
\hline
$g$ & \;$V_{g,1}/C_{g,1}\;$ & $\;V_{g,2}/C_{g,2}\;$ & $\;V_{g,3}/C_{g,3}\;$ & $\;V_{g,4}/C_{g,4}\;$\\
\hline \hline
40 & 1.019018429 & 1.002495904 & 0.9962430037 & 0.9999695265 \\
\hline
41 & 1.018547428 & 1.002435270 & 0.9963349432 & 0.9999703519 \\
\hline
42 & 1.018099193 & 1.002377513 & 0.9964224911 & 0.9999711349 \\
\hline
43 & 1.017672110 & 1.002322431 & 0.9965059531 & 0.9999718768 \\
\hline
44 & 1.017264718 & 1.002269843 & 0.9965856093 & 0.9999725812 \\
\hline
45 & 1.016875685 & 1.002219584 & 0.9966617148 & 0.9999732510 \\
\hline
46 & 1.016503797 & 1.002171501 & 0.9967345009 & 0.9999738890 \\
\hline
47 & 1.016147948 & 1.002125457 & 0.9968041809 & 0.9999744975 \\
\hline
48 & 1.015807118 & 1.002081325 & 0.9968709494 & 0.9999750781 \\
\hline
49 & 1.015480379 & 1.002038988 & 0.9969349845 & 0.9999756329 \\
\hline
\end{tabular}
\end{center}
\medskip

We see that the ratio $V_{g,n}/C_{g,n}$ apparently tends to 1 as $g\rightarrow\infty$ for any $n=1,2,3,4$
(a standard extrapolation gives 1 up to at least 6 decimal digits). For other values of $n$ the 
situation is the same. It is worth mentioning that the case $n=0$ is computationally harder
because of an additional non-trivial integration \cite{Z}, so currently we are able to 
compute $V_{g,0}$ only up to $g=30$. However, Conjecture 1 is rather accurate even in this case:

\bigskip
\noindent
\begin{center}
\begin{tabular}{|c|c|}
\hline
$g$ & \;$V_{g,0}/C_{g,0}\;$\\
\hline \hline
21 & 1.091195176 \\
\hline
22 & 1.086790774 \\
\hline
23 & 1.082792056 \\
\hline
24 & 1.079145459 \\
\hline
25 & 1.075806445 \\
\hline
26 & 1.072737684 \\
\hline
27 & 1.069907649 \\
\hline
28 & 1.067289535 \\
\hline
29 & 1.064860399 \\
\hline
30 & 1.062600498 \\
\hline
\end{tabular}
\end{center}
\medskip

Clearly, this sequence converges and its evaluated limit is 1 as well. 

Considerations that led us to Conjecture 2 are very similar to the ones described in this
section.

\subsection*{Acknowledgments} Part of this work was done at the 
Max-Planck-Institut f\"{u}r Mathematik (Bonn) in winter 2007-08, 
whose support is gratefully acknowledged. 
Special thanks are to M.~Kazarian, M.~Mirzakhani 
and D.~Zagier for numerous helpful discussions.

\end{document}